\documentclass[draft]{amsart}

\numberwithin{equation}{section}
\newcommand{\R}{{\mathbb R}}

\newcommand{\Z}{{\mathbb Z}}

\newcommand{\Hom}{\operatorname{Hom}}
\newcommand{\diag}{\operatorname{diag}}

\newcommand{\pt}{\{\text{pt}\}}

\newcommand{\nr}{N_{\R}}
\newcommand{\mr}{M_{\R}}
\newcommand{\sg}{\sigma}
\newcommand{\fan}{\Delta}

\newtheorem{thm}{Theorem}[section]
\newtheorem{cor}[thm]{Corollary}
\newtheorem{lem}[thm]{Lemma}
\newtheorem{prop}[thm]{Proposition}

\theoremstyle{remark}
\newtheorem{rem}{Remark}[section]
\newtheorem{defin}[rem]{Definition}
\begin{document}

\title
[On toric varieties and algebraic semigroups]
{On toric varieties and algebraic semigroups}

\author
{Dmitriy Boyarchenko}

\begin{abstract}
The main result of this paper is that every (separated) toric variety which 
has a semigroup structure compatible with multiplication on the 
underlying torus is necessarily affine. In the course of proving this
statement, we also give a proof of the fact that every separated toric
variety may be constructed from a certain fan in a Euclidean space. To our
best knowledge, this proof differs essentially from the ones which can
be found in the literature.
\end{abstract}

\maketitle

\section{Introduction}

In this paper we discuss two questions arising in the theory of toric
varieties and their relation to the theory of algebraic semigroups.
A {\it toric variety} is a normal algebraic variety which contains a torus as
an open dense subset, such that the natural action of the torus on itself (by
left multiplication) extends to an action of the torus on the variety.
This definition is not very enlightening, so one needs a different
description of a toric variety. It turns out that the geometry of toric
varieties is closely related to such simple objects as polyhedral cones in
Euclidean spaces. More precisely, given a fan in a Euclidean space, which is a
collection of "strongly convex rational polyhedral cones" satisfying a certain
condition, there exists a natural way of
associating a toric variety to this fan (for details, see section 2.3). In
principle,  all questions concerning the resulting toric variety may be
reformulated in terms of simple geometric statements about the fan we have
started with. This makes the study of toric varieties much easier than the
study of arbitrary algebraic varieties. However, to be able to use this
theory, we must solve the reverse problem:  does there
always exist a fan such that the associated toric variety is isomorphic to the
given one?

\medbreak

In this paper, we will see that the answer is 
affirmative in the most important case when the toric variety is {\it
separated}\footnote{The author is thankful to Professor S. Shatz for pointing
out that there exist nonseparated toric varieties.}. To justify the use of this
assumption, we will show that a toric variety associated to a fan is
necessarily separated, and we will give an example of a toric variety which is
not separated (and hence cannot be constructed from a fan). Of course, this
result is far from being original. For example, the book \cite{TV} gives a
procedure\footnote{but without checking the details} of recovering a fan from a
toric variety. However, Fulton assumes that the given toric variety {\it can}
be constructed from a fan, and then shows how the fan can be reconstructed. In
general, it seems that in most of the literature on toric varieties, the
authors try to avoid giving all the details of this reconstruction. In the
second part of this paper, we will give a {\it complete} proof of the fact
that a separated toric variety can always be constructed from a
fan. For the only
nontrivial result that we will state without a proof, the reader is referred
to the paper \cite{EC}.

\medbreak

It appears that our method differs essentially from the ones that may be found
in the literature. There are many fewer technical details to check, and our
argument makes it clear where the normality and separatedness assumptions are
used.

\medbreak

We also address a question which connects the theories of toric varieties and
algebraic semigroups. It is known that every affine toric variety has a
natural structure of an algebraic semigroup such that the multiplication
extends the multiplication on the torus.\footnote{We show this in Lemma 3.8.} A
natural question arises: is it true that every (separated) toric variety which
is an algebraic semigroup whose multiplication is compatible with that on the
torus is automatically affine?

\medbreak

The first part of the paper is devoted to proving that the answer is
affirmative in the case when the "toric semigroup" can be constructed from a
fan. Combining this result with the second part, we see that every separated
toric semigroup\footnote{For a formal definition, see 3.1.} is affine. Thus,
the study of separated toric semigroups reduces to the study of affine toric
semigroups, or, in the usual terminology, {\it toric monoids}. For a complete
description of toric monoids, the reader is referred to the paper \cite{TM}.

\medbreak

The exposition is kept as self-contained as possible.
All the relevant definitions from the theories of toric varieties and
algebraic semigroups are given. We use several basic results from the
general theory of algebraic varieties, all of which may be found in \cite{AV}.
We also use a few results on cones in Euclidean spaces. For the
proofs the reader is referred to the book \cite{TV}.

\medbreak

\section{Basic Definitions}

\subsection{Algebraic semigroups and toric varieties.}
In this paper, $k$ denotes an arbitrary algebraically closed field. We use
the term {\it algebraic variety} in the sense of \cite{AV}. In particular, we
do not assume that varieties are irreducible or separated. If $X$ is a variety
over $k$, $k[X]$ denotes its coordinate ring.

\begin{defin}
An {\it algebraic semigroup} is an algebraic variety $S$ together with a
morphism of algebraic varieties
\[
\mu: S\times S \rightarrow S,
\]
denoted multiplicatively:
\[
(a,b)\mapsto ab,
\]
such that the usual associativity condition holds:
for all $a,b,c\in S$,
\[(ab)c=a(bc).\]
A {\it zero} of an algebraic semigroup $S$ is such an element $0\in S$ (if it
exists) that for all $s\in S$, $0\cdot s=s\cdot 0=0$. If $S$ has a unit,
$G(S)$ will denote the subgroup consisting of all invertible elements of $S$. 
\end{defin}

\begin{defin}
A {\it toric variety} is a pair $(X,T)$ where $X$ is a normal algebraic variety
and $T\subset X$ is an open dense subset isomorphic to an (algebraic) torus,
such that the multiplication $T\times T\rightarrow T$ extends to an action
$T\times X\rightarrow X$ of $T$ on $X$.
\end{defin}

\begin{rem}
In this paper we are interested in algebraic semigroups containing a torus as
an open dense subset, such that the multiplication on the semigroup is
compatible with that on the torus. Let $S\supset T$ be such a semigroup. Since
the condition $s_1 s_2=s_2 s_1$ is closed and is satisfied for all $s_1,s_2\in
T$, it must be satisfied for all $s_1,s_2\in S$. Similarly, if $1$ denotes the
unit element of $T$, then the condition $1\cdot s=s\cdot 1=s$, satisfied
for all $s\in T$, must therefore be satisfied for all $s\in S$. Thus, $S$ is an
irreducible commutative algebraic semigroup with a unit. Consequently, all
elements of $T$ must be invertible. Now assume that $S$ is affine. We claim
that, in fact, $G(S)=T$. To see this, note that for an affine semigroup $S$,
there exists a closed embedding $S\hookrightarrow M_n(k)$ which is a
homomorphism of semigroups\footnote{As usual, $M_n(k)$ denotes the set of
$n\times n$ matrices with entries in $k$, and we make $M_n(k)$ into a
semigroup under {\it multiplication.}}. Moreover, if $S$ has a unit,
we may assume that it maps to the identity matrix under this embedding 
(for details, see \cite{LAM}). Hence $G(S)$ maps into $GL(n,k)$. In
particular, $G(S)$ is itself an algebraic group (the map $g\mapsto g^{-1}$ is
a morphism of algebraic varieties $G\rightarrow G$), and $G(S)$ is open in
$S$. In our situation, this means that $G(S)$ contains $T$ as an open dense
subgroup. Thus, $T=G(S)$ (see \cite{LAG}).
\end{rem}

\subsection{Fans in euclidean spaces.} Let $N$ be a lattice of finite rank
$n$, i.e. $N\equiv\Z^n$. We denote by $N_\R$ the Euclidean space $N\otimes_\Z
\R$. Let $M=\Hom(N,\Z)$ be the dual lattice and $M_\R=M\otimes_\Z \R$ the
corresponding dual space. Let $\sigma\subset\nr$ be a cone (a subset closed
under multiplication by nonnegative real numbers). We say that $\sg$ is
{\it strongly convex} if it contains no real line, {\it polyhedral} if it is an
intersection of finitely many closed half-spaces, and {\it rational} if there
exist finitely many elements $v_1,\ldots,v_k\in N$ such that 
\[
\sg=\left\{\sum_{i=1}^k \lambda_i v_i\ \text{\Huge $\mid$}\ \lambda_i\geq 0
\right\}. 
\]
We will denote by $(\cdot,\cdot):M\times N \rightarrow \Z$ and
$(\cdot,\cdot):\mr\times \nr \rightarrow \R$ the canonical pairings. Given a
cone $\sg\subset\nr$, we may form the {\it dual cone} $\sg^*:=\{u\in\mr\ \mid\ 
(u,v)\geq 0\ {\rm for\ all}\ v\in\sg\}$. If $u\in\mr$, we write
$u^{\perp}=\{v\in\nr\mid (u,v)=0\}$. Then a {\it face} of $\sg$ is a subset
which may be written as $\tau=\sg\cap u^{\perp}$ for some $u\in\sg^*$.

\begin{defin}
A {\it fan} in $\nr$ is a nonempty finite collection $\fan$ of strongly convex
rational polyhedral cones in $\nr$ satisfying the following two conditions:

\begin{enumerate}

\item If $\sigma\in\fan$ and $\tau$ is a face of $\sg$, then $\tau\in\fan$.

\item If $\sigma,\tau\in\fan$, then $\sigma\cap\tau$ is a face of both $\sg$
and $\tau$.

\end{enumerate}
If there exists a cone $\sg\in\fan$ such that $\fan$ consists of $\sg$
together with all of its faces, we say that the fan $\fan$ is {\it generated
by} the cone $\sg$.

\end{defin}

\subsection{Fans and toric varieties.}
Now we will show how to associate, to each fan $\fan$ in $\nr$, a toric
variety. First, for each $\sg\in\fan$, we define an affine variety $U_\sg$.
Let $\sg^*$ be the dual cone and let $S_\sg$ be the subsemigroup of $M$
defined by $S_\sg=\sg^*\cap M$. We may form the {\it semigroup ring}
$k[S_\sg]$\footnote{It is the commutative $k$-algebra generated by the
"monomials" $\{\chi^v\mid v\in S_\sg\}$ subject to the relation
$\chi^v\cdot\chi^w=\chi^{v+w}$. In the following discussion, we identify
elements of $S_\sg$ with the corresponding monomials in $k[S_\sg]$.}. Clearly,
it is a finitely generated $k$-algebra without zero divisors, so we set
$U_\sg=\operatorname{Spec}(k[S_\sg])$ (technically, we should write m-Spec
here, since we consider only the closed points of the corresponding scheme).
Thus $U_\sg$ is an irreducible algebraic variety over $k$. 

Next, we patch the affines $U_\sg$ together to obtain a toric variety. For all
$\sg,\tau\in\fan$, the inclusion maps $\sg\cap\tau\hookrightarrow\sg$ and
$\sg\cap\tau\hookrightarrow\tau$ correspond to maps
$\sg^*\hookrightarrow(\sg\cap\tau)^*$, $\tau^*\hookrightarrow(\sg\cap\tau)^*$,
and hence to homomorphisms of semigroups $S_\sg\hookrightarrow
S_{\sg\cap\tau}$, $S_\tau\hookrightarrow
S_{\sg\cap\tau}$. These induce homomorphisms of $k$-algebras
$k[S_\sg]\hookrightarrow
k[S_{\sg\cap\tau}]$, $k[S_\tau]\hookrightarrow
k[S_{\sg\cap\tau}]$, and, finally, dominant morphisms of algebraic varieties
\[
U_{\sg\cap\tau}\rightarrow U_\sg
\]
and
\[
U_{\sg\cap\tau}\rightarrow U_\tau.
\]
It turns out that these morphisms are open embeddings, so we may glue the
affines $U_\sg$, $\sg\in\fan$, along their open subvarieties
$U_{\sg\cap\tau}$, to obtain an irreducible variety $X$. Note that $\{0\}$ is a
face of every strongly rational cone, so by condition 1 of definition 2.4,
$\{0\}\in\fan$. Moreover, for each $\sg\in\fan$, the inclusion
$\{0\}\hookrightarrow\sg$ induces an open embedding \[
U_{\{0\}}\hookrightarrow U_\sg.
\]
But $S_{\{0\}}=M$, hence we have natural isomorphisms 
\[
k[S_{\{0\}}]\cong
k[T_1,T_1^{-1},\ldots,T_n,T_n^{-1}]
\]
and
\[
U_{\{0\}}\cong T=(k^*)^n.
\]
Thus all the affines $U_\sg$ contain the torus $T$ as an open dense subset,
and hence so does $X$. It remains to define the action of $T$ on $X$. We first
construct this action on each of the $U_\sg$'s. To show that the
multiplication morphism $T\times T\rightarrow T$ extends to a morphism
$T\times U_\sg\rightarrow U_\sg$, it suffices to prove the corresponding
statement on the level of coordinate rings, namely, that we have the
commutative diagram
\[
\begin{array}{ccc}
k[T] & \longrightarrow & k[T]\otimes_k k[T] \\
\cup & & \cup \\
k[U_\sg] & \longrightarrow & k[T]\otimes_k k[U_\sg] \\
\end{array}
\]
Note that the lattice $M$ is naturally isomorphic to the lattice of characters
of $T$. Hence $k[U_\sg]=k[S_\sg]$ has a basis consisting of characters of $T$.
But if $\chi$ is a character of $T$, then the homomorphism
\[
k[T]\rightarrow k[T]\otimes_k k[T]
\]
maps $\chi$ to $\chi\otimes\chi$. This argument shows that, in fact, we have a
diagram
\[
\begin{array}{ccc}
k[T] & \longrightarrow & k[T]\otimes_k k[T] \\
\cup & & \cup \\
k[U_\sg] & \longrightarrow & k[U_\sg]\otimes_k k[U_\sg] \\
\end{array}
\]
Since the homomorphism in the top row of this diagram is coassociative, so is
the homomorphism in the bottom row. Thus, $U_\sg$ has a natural semigroup
structure which is compatible with the multiplication on the torus.

It is easy to check that the actions of the torus on the different $U_\sg$'s
agree, and hence give an action of $T$ on $X$. This completes the construction
of the toric variety $X$. For more details, see \cite{TV}.

\begin{defin}
Given a finite dimensional lattice $N$ and a fan $\fan$ in $\nr=N\otimes_\Z
\R$, the toric variety $X$ constructed above will be called the {\it toric
variety associated with the fan} $\fan$, written $X=X(\fan)$.
\end{defin}

Observe that if the fan $\fan$ is generated by a single cone $\sigma\in\fan$,
then the associated toric variety is affine. Later we will prove the converse
of this statement (Corollary 3.9).

\section{Main Results}

In this section, we state and prove our main results and their
corollaries. 

\subsection{Toric semigroups.}
First, we define the main object of our study.
\begin{defin}
A {\it toric semigroup} is a toric variety $X\supset T$ which has a structure
of an algebraic semigroup compatible with the multiplication on the torus.
\end{defin}
\begin{rem}
In modern literature (cf. \cite{TM}), affine toric semigroups are usually
called {\it toric monoids}.
\end{rem}

\noindent
We intend to prove

\begin{thm}
Let $N$ be a finite dimensional lattice, $\fan$ a fan in $N_\R$ and
$X=X(\fan)$ the associated toric variety. If $X$ has a structure of a toric
semigroup, then the fan $\fan$ is generated by a single cone $\sg\in\fan$. In
particular, $X$ is affine.
\end{thm}

We will use several lemmas. We call a cone
$\sg\subset\nr$ {\it nondegenerate} if its interior is nonempty, or,
equivalently, if $\nr$, as a real vector space, is spanned by $\sg$. In the
previous section, we have shown that for all $\sg\in\fan$, the corresponding
affine variety $U_\sg$ has a natural structure of a toric monoid. We have

\begin{lem}
The toric monoid $U_\sg$ has a zero if and only if the cone $\sg$ is
nondegenerate. 
\end{lem}

\noindent
{\it Proof.} See Appendix, section 4.1.

\medbreak

To use this result, we would like to be able to relate the semigroup
structure on $X$ with those on the open affines $U_\sg$. For this, we apply

\begin{lem}
The multiplication morphism $X\times X\rightarrow X$ maps $U_\sg\times U_\sg$
into $U_\sg$, and its restriction to a morphism 
\[
U_\sg\times U_\sg\rightarrow U_\sg
\]
coincides with the natural monoid structure on $U_\sg$. More generally, let
$X\supset T$ be a separated toric variety and let $Y\subset X$ be a 
$T$-invariant subvariety
containing $T$ (so $Y$ is naturally a toric variety). If both $X$ and
$Y$ have structures of toric semigroups, the two structures agree
on $Y$.
\end{lem}

\noindent
{\it Proof.} See Appendix, section 4.2.

\smallbreak

\begin{cor}
Under the assumptions of Theorem 3.1,
the fan $\fan$ contains at most one nondegenerate cone.
\end{cor}

\begin{proof} 
Suppose, to obtain a contradiction, that there are two distinct
nondegenerate cones $\sg,\tau\in\fan$. By Lemma 3.3, the semigroup structure
on $X$ induces monoid structures on $U_\sg$ and $U_\tau$, and by Lemma 3.2,
both have zeroes, say $0_\sg\in U_\sg$ and $0_\tau\in U_\tau$. The relation
$x\cdot 0_\sg=0_\sg\cdot x=0_\sg$, which holds for all $x\in T$ (because
$T\subset U_\sg$), must therefore hold for all $x\in X$. Thus, $0_\sg$ is a
zero of $X$. Similarly, $0_\tau$ is also a zero of $X$. But a semigroup cannot
have two distinct zeroes, so $0_\sg=0_\tau$. Consequently,
\[
0_\sg\in U_\sg\cap U_\tau=U_{\sg\cap\tau},
\]
and in particular $U_{\sg\cap\tau}$ has a zero. But $\sg,\tau$ are distinct
cones, whence $\sg\cap\tau$ is a proper face of $\sg$, and is
therefore degenerate. This contradicts Lemma 3.2.
\end{proof}

\medbreak

This result already excludes many possibilities for the fan $\fan$. It can be
used indirectly to exclude even more possibilities. For example, suppose that
all the cones in $\fan$ are degenerate, and in fact are contained in a certain
subspace $V\subset\nr$, where $V$ has a basis consisting of elements of $N$.
Then, clearly, $\fan$ can be considered as a fan in $V$. Hence, if $\fan$
contains two cones which are nondegenerate as cones in $V$, we may again apply
Corollary 3.4.\footnote{We will work out the details of this argument after
Lemma 3.5.} 

\medbreak

Still, observe that there do exist fans $\fan$ which are not
generated by one cone, and to which our previous argument does not apply. For
example, think of a collection of degenerate cones in $\nr$ whose union spans
$\nr$. To exclude such possibilities, we apply

\begin{lem}
Under the assumptions of Theorem 3.1, the union of all cones in $\fan$, viewed
as a subset on $\nr$, is closed under addition.
\end{lem}

\noindent
{\it Proof.} See Appendix, section 4.3.

\medbreak

Now we may finish the proof of Theorem 3.1. First, note that we may assume
that the union of all the cones in $\fan$ span the whole space $\nr$.
Indeed, let
\[
V=\operatorname{span}_{\R}\left(\cup_{\sg\in\fan} \sg\right).
\]
Since the cones in $\fan$ are rational, $V$ has a basis consisting of elements
on $N$. Thus, if we set $N_\fan=N\cap V$, then $V=N_\fan\otimes_\Z \R$, and we
may consider $\fan$ as a fan in $V$. Since $V$ is a vector space, $N_\fan$ is
a {\it saturated subsemigroup} of $N$, i.e. if $v\in N$, $m\in\Z$ and $m\cdot
v\in N_\fan$, then $v\in N_\fan$. Therefore we may choose a sublattice
$N'_\fan\subset N$ such that $N=N_\fan\oplus N'_\fan$. Corresponding to this
decomposition, we have $T=T_\fan\times T'_\fan$ where $T_\fan$ and $T'_\fan$
are the tori associated to the lattices $N_\fan$ and $N'_\fan$, respectively.
Hence, $X=Y\times T'_\fan$ where $Y$ is the toric variety associated to the
fan $\fan$ viewed as a fan in $V$. Since $T_\fan\times\{1\}$ is dense in
$Y\times\{1\}$, the multiplication morphism $X\times X\rightarrow X$ restricts
to a morphism 
\[
(Y\times\{1\})\times (Y\times\{1\})\rightarrow (Y\times\{1\}),
\]
thus making $Y$ into a toric semigroup. But now the fan $\fan$ spans the
corresponding space $V$, by construction.

\medbreak

Now write $\sg$ for the union of all cones in $\fan$. Then $\sg$ is a cone,
and it is convex because it is closed under addition (Lemma 3.5). Since $\sg$
generates $\nr$ as a real vector space, $\sg$ must have nonempty interior.
This already implies that $\fan$ contains at least one nondegenerate cone,
$\tau$, since otherwise $\sg$ is a finite union of nowhere dense sets. In
fact, since $\fan$ has at most one nondegenerate cone (Corollary 3.4), we must
have $\tau=\sigma$, otherwise $\sigma\setminus\tau$ has nonempty interior
(because $\tau$ is closed). But now for each $\tau'\in\fan$, we have
$\tau'\subset\sigma=\tau$ by construction, whence, by condition 2 of
Definition 2.2.1, $\tau'$ is a face of $\tau$. Consequently, $\fan$ is
generated by the cone $\tau$.

\hskip11.1cm
Q.E.D.

\bigbreak

\subsection{Construction of fans from toric varieties.} Here we consider the
other problem mentioned in the introduction. Given a toric variety $X\supset
T$, does there exist a fan $\fan$ such that $X\cong X(\fan)$? First let us show
that it is necessary to assume that $X$ is separated:

\begin{prop}
If $N$ is a finite dimensional lattice and $\fan$ is a fan in $\nr$, then the
associated toric variety $X(\fan)$ is separated.
\end{prop}

\noindent
{\it Proof.} See Appendix, section 4.4.

\medbreak

On the other hand, there do exist nonseparated toric varieties. The easiest
example to construct is the standard example of a nonseparated variety: the
affine line with a doubled origin (see \cite{AV}, the example after Lemma
3.3.4). We consider two copies of the affine line, ${\mathbb A}^1_1$ and
${\mathbb A}^1_2$, and glue them along their open subsets ${\mathbb
A}^1_1\setminus\{0\}$ and ${\mathbb A}^1_2\setminus\{0\}$, with the
identification morphism  \[
{\mathbb A}^1_1\setminus\{0\}\rightarrow 
{\mathbb A}^1_2\setminus\{0\}
\]
given by the identity map. The resulting nonseparated variety $X$ contains the
one dimensional torus $T={\mathbb A}^1\setminus\{0\}$, which acts by
multiplication on both copies of the affine line, and hence on $X$.

\begin{rem}
In the following discussion, we will implicitly assume the following criterion
of separability:

\smallbreak

{\it an algebraic variety $X$ is separated if and only if a morphism
$U\rightarrow X$ from 

an open dense subset $U$ of a smooth curve $C$ extends
in at most one way 

to a morphism $C\rightarrow X$ (see \cite{AV}, Proposition 7.2.2),}

\smallbreak

\noindent
more precisely, the implication $(\Rightarrow)$ for the special case
$C={\mathbb A}^1$, $U={\mathbb A}^1\setminus\{0\}$.
\end{rem}

\bigbreak

According to the discussion above, the best statement we can try to prove is

\begin{thm}
Let $X\supset T$ be a separated toric variety. Let $N=\Hom(k^*,T)$ be
the lattice of one-parameter subgroups of $T$, and $M=N^*=\Hom(T,k^*)$ the
dual lattice of characters of $T$. Then there exists a fan $\fan$ in
$\nr=N\otimes_\Z \R$ such that the natural isomorphism 
\[
\operatorname{Spec}(k[M])
\begin{array}{c}
_{\cong} \\
\longrightarrow \\
 \\
\end{array}
T
\]
extends to an isomorphism 
\[
X(\fan)
\begin{array}{c}
_{\cong} \\
\longrightarrow \\
 \\
\end{array}
X.
\]
\end{thm}

\medbreak

First, let's prove the theorem in the particular case when $X$ is affine.
The key idea is to observe that $X$ has a unique structure
of a toric monoid (NB: earlier we proved this statement only assuming that
$X=U_\sg$ for some strongly convex rational polyhedral cone $\sg$) and look
at the "one-parameter subsemigroups" of $X$. We have

\begin{lem}
Let $X\supset T$ be an affine toric variety. Then there exists a
structure of a toric monoid on $X$.
\end{lem}

\noindent
{\it Proof.} See Appendix, section 4.5.

\medbreak

\noindent
As a by-product of our discussion, we obtain a statement, which is not at all
obvious from the definitions.

\begin{cor}
Let $N$ be a finite dimensional lattice and $\fan$ a fan in $\nr$. The
associated toric variety $X=X(\fan)$ is affine if and only if $\fan$ is
generated by a single cone.
\end{cor}
\begin{proof}
The direction $(\Leftarrow)$ is clear. For $(\Rightarrow)$, combine Lemma 3.8
and Theorem 3.1.

\end{proof}

\medbreak

\noindent
Now we introduce an important

\begin{defin}
If $X$ is an algebraic semigroup, a {\it one-parameter subsemigroup} of $X$ is
a homomorphism of algebraic semigroups $k\rightarrow X$, where $k$ is viewed
as an algebraic semigroup {\it under multiplication}. The set of all such is
denoted by $\Hom(k,X)$. In case $X$ is commutative, $\Hom(k,X)$ has a natural 
structure of an abelian semigroup under pointwise multiplication.
\end{defin}

\begin{rem}
The notion of a one-parameter subsemigroup corresponds to the well-known
notion of a one-parameter subgroup of an algebraic group (see, for example,
\cite{LAG}). In fact, if $S$ is an irreducible affine algebraic semigroup, so
that the group of invertible elements of $S$, $G(S)$, is itself an algebraic
group, then the restriction of every one-parameter subsemigroup $k\rightarrow
S$ to $k^*$ is a one-parameter subgroup of $G(S)$. Moreover, a subsemigroup
$k\rightarrow S$  is determined by its restriction to $k^*$ (since affine
varieties are separated), and hence if $S$ is abelian, it is possible to view
$\Hom(k,S)$ as a subsemigroup of $\Hom(k^*,G(S))$.
\end{rem}

\medbreak

With these tools, we may complete the proof of Theorem 3.7 in the affine case.
In fact, we prove a slightly more precise statement.

\begin{prop}
Under the conditions of Theorem 3.7, assume that $X$ is affine, and view $X$
as a toric monoid (cf. Lemma 3.8). Let $\sg$ be the cone
generated\footnote{That is, $\sg$ is the convex hull of $\Hom(k,X)$.} by the
subsemigroup $\Hom(k,X)$ of the lattice $N=\Hom(k^*,T)$ (cf. Remark 3.2.1) and
let $\fan$ be the fan in $\nr$ generated by the cone $\sg$ (see Definition
2.2.1). Then the conclusion of Theorem 3.7 is valid for the fan $\fan$.

\end{prop}

\noindent
{\it Proof.} See Appendix, section 4.6.

\bigbreak

\noindent
Now we attack the general case. The following result is very useful for us

\begin{prop}
Let $T$ be a torus and let $X$ be a normal variety on which $T$ acts. Then,
for any point $x\in X$, there is a $T$-invariant affine open neighborhood of
$x$.
\end{prop}
\begin{proof}
See \cite{EC}, Corollary 2 of Lemma 8.
\end{proof}

\begin{cor}
Under the assumptions of Theorem 3.7, there exists a finite open cover of $X$
by affine toric varieties $V_i\supset T$.
\end{cor}
\begin{proof}
By Proposition 3.11, there exists a finite cover of $X$ by invariant
open affine subsets. But every open subset of $X$ intersects $T$, since $T$ is
dense in $X$, and an invariant subset of $X$ intersecting $T$ must contain $T$,
because $T$ itself is an orbit under the action of $T$.
\end{proof}

\medbreak

Now we fix a finite open cover $\{V_i\}$ of $X$ by affine toric varieties as
in Corollary 3.12. By Proposition 3.10, for each of the $V_i$ we have the
corresponding strongly convex rational polyhedral cone $\sg_i\subset\nr$. Let
$\fan$ be the finite set of cones consisting of the cones $\sg_i$ together
with all their faces. We claim that $\fan$ is a fan and that the conclusion of
Theorem 3.7 is valid for $\fan$.

\begin{rem}
One needs to be a little careful here. In this construction, we implicitly use
the assumption that $X$ is separated. The reason is that it makes sense to
consider the cones $\sg_i$ simultaneously only if we know that every
one-parameter subgroup of $T$ extends in {\it at most one way} to a
one-parameter subsemigroup of some $V_i$. To illustrate this phenomenon,
consider again the affine line with a doubled origin (see the example after
Proposition 3.6). It has an open cover by affine toric varieties, namely, the
two copies of the affine line that were used for gluing. But the cones
corresponding to these two affine varieties coincide! So in the nonseparated
case, intersection of cones does not necessarily correspond to intersection of
open affine varieties. In a moment we'll see that this phenomenon never occurs
in the separated case.
\end{rem}

\medbreak

To show that $\fan$ is a fan, it will suffice to see that for all $i,j$,
$\sg_i\cap\sg_j$ is a face of both $\sg_i$ and $\sg_j$. Consider
$V_{ij}=V_i\cap V_j$. It is certainly an open subset of $X$ containing $T$.
But in fact, since $X$ is separated, $V_{ij}$ is affine, and the natural
$k$-algebra homomorphism $k[V_i]\otimes_k k[V_j]\rightarrow k[V_{ij}]$ is
surjective (see \cite{AV}, Proposition 3.3.5). We know that
$k[V_i]=k[S_{\sg_i}]$ where we use the same notation as in section 2.3. We
also have the cone $\sg_{ij}$ corresponding to the affine toric variety
$V_{ij}$ (Proposition 3.10). By separability, a one-parameter subgroup
of $T$ extends to a one-parameter subsemigroup of $V_{ij}$ if and only if it
extends to one-parameter subsemigroups both of $V_i$ and of
$V_j$\footnote{Here we implicitly use the fact that the structures of toric
semigroups on $V_i$ and $V_j$ agree on $V_{ij}$. In fact, this follows from
the second statement of Lemma 3.3}. That is,  \[
\Hom(k,V_{ij})=\Hom(k,V_i)\cap\Hom(k,V_j), \] as subsemigroups of
$\Hom(k^*,T)$. Consequently, $\sg_{ij}=\sg_i\cap\sg_j$. Say we already know
that $\sg_{ij}$ is a face of both $\sg_i$ and $\sg_j$. Then $\fan$ is a fan,
and we are done, because the toric variety $X(\fan)$ associated to $\fan$ is
obtained by gluing the affine varieties corresponding to the cones $\sg_i$
along their open subvarieties corresponding to the cones $\sg_{ij}$, and the
induced morphism $X(\fan)\rightarrow X$ (which exists by the universal
property of gluing) is clearly an isomorphism.

\medbreak

To see that $\sg_{ij}$ must be a face of $\sg_i$ (and hence also $\sg_j$, by
symmetry), we make the following observation. In any case, the inclusion
$\sg_{ij}\subset\sg_i$ corresponds to the inclusion $\sg_i^*\subset\sg_{ij}^*$,
hence to the inclusion $S_{\sg_i}\subset S_{\sg_{ij}}$, which induces a
dominant morphism $U_{\sg_{ij}}\rightarrow U_{\sg_i}$. We know that if
$\sg_{ij}$ {\it is} a face of $\sg_i$, then this morphism is an open embedding.
We conclude the proof with

\begin{lem}
If $\sg_{ij}$ is not a face of $\sg_{i}$, then the morphism
$U_{\sg_{ij}}\rightarrow U_{\sg_i}$ is not an open embedding.
\end{lem}

\medbreak

Note that this observation automatically implies that $\sg_{ij}$ is a face of
$\sg_i$, because in our case the morphism $U_{\sg_{ij}}\rightarrow U_{\sg_i}$
is just the inclusion $V_{ij}\hookrightarrow V_i$. For the proof of the lemma,
see Appendix, section 4.7.

\section{Appendix}

\subsection{Proof of Lemma 3.2.}
The key idea here is to write down the definition of a zero of an affine
semigroup in terms of homomorphisms of coordinate rings. Since $U_\sg$ is
commutative, an element of $U_\sg$ is a zero if an only if it is a left zero.
Note that an element of $U_\sg$ can be naturally viewed as a morphism
$\pt\rightarrow U_\sg$ where $\pt=\operatorname{Spec} k$ is the one-point
variety. Therefore we may restate our definition of zero as follows: it is a
morphism $0:\pt\hookrightarrow U_\sg$ such that the diagram
\[
\begin{array}{ccccc}
\pt\times U_\sg & \longrightarrow & U_\sg\times U_\sg & \longrightarrow & U_\sg
\\
\parallel & & & & \parallel \\
\pt\times U_\sg & \longrightarrow & \pt & \longrightarrow & U_\sg \\
\end{array}
\]
commutes, where the arrows are (from the left to the right) in the top row:
$0\times\operatorname{id}_{U_\sg}$, multiplication; in the bottom row:
projection onto the factor $\pt$, $0$. In terms of coordinate rings, a zero of
$U_\sg$ is such a homomorphism of $k$-algebras $f:k[U_\sg]\rightarrow k$ that
the following diagram commutes:
\[
\begin{array}{ccccc}
k\otimes_k k[U_\sg] & \longleftarrow & k[U_\sg]\otimes_k k[U_\sg] &
\longleftarrow & k[U_\sg] \\
\parallel & & & & \parallel \\
k\otimes_k k[U_\sg] & \longleftarrow & k & \longleftarrow & k[U_\sg] \\
\end{array}
\]
and now the arrows are (from the left to the right) in the top row:
$f\otimes\operatorname{id}_{k[U_\sg]}$, comultiplication; in the bottom row:
natural inclusion of $k$ into a $k$-algebra, $f$.

\medbreak

We recall that $k[U_\sg]=k[S_\sg]$. I claim that the diagram above exists if
and only if there exists a $k$-algebra homomorphism $f:k[S_\sg]\rightarrow k$
which kills all nonzero elements of $S_\sg$. Indeed, let $v\in S_\sg$, then in
the top row it maps first to $v\otimes v$, and then to $f(v)\otimes v$. In the
bottom row, it maps first to $f(v)$, and then to $f(v)\otimes 1$. Clearly,
$f(v)\otimes v=f(v)\otimes 1$ in $k\otimes_k k[S_\sg]$ if and only if either
$v=0$, $f(v)=1$, or $f(v)=0$. So if the diagram above commutes, then $f$ has
the required property. Conversely, if $f$ has the required property, then the
diagram commutes, because $S_\sg$ is a basis of $k[S_\sg]$.

\medbreak

On the other hand, there exists a homomorphism $f$ with the required property
if and only if there exists a map $g:S_\sg\rightarrow k$ such that $g(0)=1$,
$g(v)=0$ for all $v\neq 0$, and $g(v+w)=g(v)g(w)$ for all $v,w\in S_\sg$.
If $\sg$ is degenerate, then there are nonzero elements $v,w\in S_\sg$ such
that $v+w=0$, and then such a map $g$ cannot exist (applying $g$ to the
equality $v+w=0$ gives $0\cdot 0=1$, an absurd). On the other hand, if $\sg$
is nondegenerate, then such elements $v,w$ do not exist, so if we define $g$
by $g(0)=1$, $g(v)=0$ for $v\neq 0$, then $g$ satisfies $g(v+w)=g(v)g(w)$ for
all $v$, $w$. This completes the proof.

\medbreak

\subsection{Proof of Lemma 3.3.}
By Proposition 3.6 (whose proof does not depend on this lemma), $X$
is separated, so the first statement is a special case of the second. 
We have
two morphisms 
\[
\phi_1,\phi_2 : Y\times Y\rightarrow X,
\]
namely, $\phi_1$ is the composition
\[
Y\times Y\hookrightarrow X\times X\rightarrow X,
\]
the second morphisms being the multiplication on $X$, and $\phi_2$ is the
composition 
\[
Y\times Y\rightarrow Y\hookrightarrow X,
\]
the first morphism being the multiplication on $Y$ 
2.3. Since the multiplications on both $X$ and $Y$ agree with the
multiplication on the torus, $\phi_1$ and $\phi_2$ agree on $T\times T$. But
$X$ is separated and $T\times T$ is dense in 
$Y\times Y$, which implies that $\phi_1=\phi_2$. Thus, the diagram
\[
\begin{array}{ccc}
Y\times Y & \longrightarrow & Y \\
    \cap  &       &   \cap \\
X\times X & \longrightarrow & X \\
\end{array}
\]
commutes, which is exactly the statement of Lemma 3.3.

\medbreak

\subsection{Proof of Lemma 3.5.}
Choose $\sg\in\fan$ and a point
$v\in N$. Note that since we constructed the torus $T$ as $\operatorname{Spec}
k[M]$, the lattice $M$ is naturally identified with the character lattice of
$T$, and then $N$ is identified with the lattice of one-parameter subgroups of
$T$. Thus, $v$ naturally corresponds to a one-parameter subgroup, $s_v$, of
$T$. Let us first check that $v\in\sg$ if and only if $s_v$ extends to a
one-parameter subsemigroup of $U_\sg$ (see Definition 3.2.1). On the level of
coordinate rings, $s_v$ extends to a one-parameter subsemigroup of $U_\sg$ if
and only if under the homomorphism $s_v^*:k[T]\rightarrow k[k^*]=k[x,x^{-1}]$,
the subring $k[U_\sg]$ of $k[T]$ is mapped into $k[x]$. We also note that on
elements of $M$ (naturally viewed as a subset of $k[T]=k[M]$), the homomorphism
$s_v^*$ is given by $\chi\mapsto x^{(\chi,v)}$. Thus $s_v$ extends to a
one-parameter subsemigroup of $U_\sg$ if and only if $v$ has nonnegative
pairing with all elements of $S_\sg$. By definition, this happens if and only
if $v\in(\sg^*)^*$. But $\sg=(\sg^*)^*$ (see \cite{TV}).

\medbreak

Now it's not hard to complete the proof. We know that if elements $v,w\in N$
correspond to one-parameter subgroups $s_v,s_w$ of $T$, then $v+w$ corresponds
to $s_v\cdot s_w$ (pointwise product of the subgroups). Now if $v\in\sg$,
$w\in\tau$ for some cones $\sg,\tau\in\fan$, we know that both $s_v$ and $s_w$
extend to one-parameter subsemigroups of $U_\sg$ and $U_\tau$, respectively,
and hence to subsemigroups of $X$, because the multiplication on $X$ is
compatible with those on $U_\sg$ and $U_\tau$ (see Lemma 3.3). We denote the
extensions again by $s_v$ and $s_w$. Then, clearly, $s_{v+w}$ also extends to
a subsemigroup of $X$, namely, $s_v\cdot s_w$. There is a cone
$\psi\in\fan$ such that $(s_v\cdot s_w)(0)\in U_\psi$, and since
$s_{v+w}(k^*)\subset T$, the image of the extension, $s_v\cdot s_w$, of
$s_{v+w}$ is contained in $U_{\psi}$. Thus, $v+w\in\psi$, completing the proof.

\medbreak

\subsection{Proof of Proposition 3.6.}
Write $\diag : X\hookrightarrow X\times X$ for the diagonal embedding. We must
prove that $\diag(X)$ is closed in $X\times X$. By construction, $X=X(\fan)$ is
covered by the affines $U_\sg$. Hence, it is enough to show that for all
$\sg,\tau\in\fan$, $\diag(U_\sg\cap U_\tau)$ is closed in $U_\sg\times
U_\tau$. We know that $U_\sg\cap U_\tau=U_{\sg\cap\tau}$ is affine, so it
suffices to see that the induced $k$-algebra homomorphism
\[
\diag^* : k[U_\sg]\otimes_k k[U_\tau] \rightarrow k[U_{\sg\cap\tau}]
\]
is surjective. But, in the notation of section 2.3, $k[U_\sg]=k[S_\sg]$ and
$k[U_\tau]=k[S_\tau]$, and since addition in the lattice $M$ corresponds to
multiplication in the algebras $k[S_\sg]$, $k[S_\tau]$, the image of $\diag^*$
is $k[S_\sg+S_\tau]$. Since $(\sg\cap\tau)^*=\sg^* + \tau^*$, the proof is
complete.

\medbreak

\subsection{Proof of Lemma 3.8.}
It suffices to prove that there exists a commutative diagram
\[
\begin{array}{ccc}
k[T] & \rightarrow & k[T]\otimes_k k[T] \\
\cup & & \cup \\
k[X] & \rightarrow & k[X]\otimes_k k[X] \\
\end{array}
\]
because this will give a morphism $X\times X\rightarrow X$, compatible with
the multiplication on the torus, and associativity will follow from the
coassociativity of the bottom row of the diagram, which in turn is implied by
the coassociativity of the top row. We know, by definition, that there exists
a commutative diagram of the required form if we replace $k[X]\otimes_k k[X]$
by $k[T]\otimes_k k[X]$. But the torus is commutative, so the left action of
$T$ on $X$ induces a right action of $T$ on $X$, which is also compatible with
the multiplication on the torus. Hence, in the diagram above, $k[X]\otimes_k k[X]$
may also be replaced by $k[X]\otimes_k k[T]$. But, by simple linear algebra,
\[
k[X]\otimes_k k[T]\cap k[T]\otimes_k k[X]=k[X]\otimes_k k[X],
\]
as subalgebras of $k[T]\otimes_k k[T]$. This completes the proof.

\subsection{Proof of Proposition 3.10.}
The left action of $T$ on itself (by multiplication) naturally induces an
action of $T$ on the coordinate ring $k[T]$. The latter decomposes into a
direct sum of one-dimensional $T$-invariant subspaces, namely, those spanned by
the characters of $T$. Since the action of $T$ on itself extends to an action
of $T$ on $X$, the subalgebra $k[X]$ of $k[T]$ is $T$-invariant, and hence
itself is a direct sum of one-dimensional subspaces spanned by characters of
$T$. Let $S$ denote the set of all characters of $T$ contained in $k[X]$, so
that $k[X]=k[S]$. Since $k[X]$ is closed under multiplication, $S$ is closed
under addition. Thus, if we let $\sg$ be the cone in $\mr$ generated by $S$
and $\tau=\sg^*$, then $\sg=\tau^*$ and therefore $X=U_\tau$. We already know
(see the proof of Lemma 3.5) that $\sg$ is generated by the set of
one-parameter subsemigroups of $U_\sg$. This completes the proof.

\medbreak

\subsection{Proof of Lemma 3.13.} To simplify notation, write $\sg$ for
$\sg_i$ and $\tau$ for $\sg_{ij}$.  We recall that, in general, given an
affine variety $V$, every open affine subset of $V$ has the form $\{x\in V\mid
f(x)\neq 0\}$ for some $f\in k[V]$. For us, this says that $U_\tau\rightarrow
U_\sg$ is an open embedding if and only if the corresponding homomorphism of
coordinate rings $k[S_\sg]\hookrightarrow k[S_\tau]$ may be realised as an
embedding $k[S_\sg]\hookrightarrow k[S_\sg]_{(f)}$ for some $f\in k[S_\sg]$,
that is (if we now view $k[S_\sg]$ as a subalgebra $k[S_\tau]$), if and only
if 
\[
k[S_\tau]=k[S_\sg,f^{-1}]
\]
for some $f\in k[S_\sg]$. Now we assume that such $f$ exists. Then, in
particular, $f$ is nonzero on $T$, and hence, multiplying $f$ by a nozero
constant, we may assume that $f$ is a character of $T$. Consequently, we have
$f\in S_\sg$. Now $S_\tau=S_\sg+\Z f$. Therefore, 
\[
\tau=\{v\in\sg\mid (-f,v)\geq 0\}=\{v\in\sg\mid (f,v)=0\}
\]
is a face of $\sg$, completing the proof.

\end{document}